\documentclass[12pt]{article}

\usepackage{amsthm,amsmath,amssymb}

\usepackage{graphicx}
 \usepackage{epsfig}
  \usepackage{latexsym, epsfig, psfrag,eepic,colordvi,bm}
  \usepackage{graphics,color}
  \usepackage{graphics}
  \usepackage{graphicx}
  \usepackage{epsfig}
  \usepackage[all]{xy}
  
  \usepackage{amssymb}
  \usepackage{amsthm,amsmath}
  
  \usepackage{latexsym, epsfig, psfrag,eepic,colordvi,bm}
  \usepackage{graphics,color}
  \usepackage{amsmath,amsfonts,amssymb,amscd}
  \usepackage[all]{xy}
  \usepackage{epstopdf}
  
  \usepackage{amsthm,amsmath,amssymb}
  
  \usepackage{graphicx}
  \usepackage{amsthm,amsmath,amssymb}
  
  \usepackage{graphicx,epsfig}
\usepackage[colorlinks=true,citecolor=black,linkcolor=black,urlcolor=blue]{hyperref}

\usepackage[all]{xy}
\usepackage{float}

\theoremstyle{plain}
\newtheorem{theorem}{Theorem}[section]

\newtheorem{lemma}{Lemma}[section]
\newtheorem{corollary}{Corollary}[section]
\newtheorem{proposition}{Proposition}[section]

\theoremstyle{definition}

\theoremstyle{remark}
\newtheorem{remark}[theorem]{Remark}

\newcommand{\se}{\textrm{E}}

\date{}

\title{\bf New formulas counting one-face maps and Chapuy's recursion}

\author{ Ricky X. F. Chen$^1$, Christian M. Reidys$^2$\\
	\small  Biocomplexity Institute and Dept. of Mathematics,\\[-0.8ex]
	\small Virginia Tech, 1015 Life Science Circle,\\[-0.8ex]
	\small Blacksburg, VA 24061, USA\\
	\small\tt $^1$cxiaof6@vt.edu, $^2$duck@santafe.edu
}

\begin{document}

\maketitle

\begin{abstract}
	In this paper, we begin with the Lehman-Walsh formula counting one-face maps
	and construct two involutions
	on pairs of permutations to obtain a new formula for the
	number $A(n,g)$ of one-face maps
	of genus $g$. Our new formula is in the form of a convolution of the Stirling numbers 
	of the first kind which immediately implies a formula for the generating function
	$A_n(x)=\sum_{g\geq 0}A(n,g)x^{n+1-2g}$ other than the well-known Harer-Zagier formula.
	By reformulating our expression 
	for $A_n(x)$ in terms of the backward shift 
	operator $\se: f(x)\rightarrow f(x-1)$ and proving a property satisfied by polynomials of the form
	$p(\se)f(x)$, we easily establish the recursion obtained by Chapuy for
	$A(n,g)$. Moreover, we give a simple combinatorial interpretation for the Harer-Zagier recurrence.



\end{abstract}



\section{Introduction}
A one-face map is a graph embedded in a closed orientable surface such
that the complement is homeomorphic to an open disk. It is well known that the combinatorial counterpart
of one-face maps are fatgraphs with one boundary component~\cite{edmonds}.
A fatgraph is a graph with a specified cyclic order of the ends of edges incident to each vertex of the graph. 
In this paper we will not differentiate between maps and fatgraphs.

\begin{figure}[!htb]
	\centering
	\includegraphics[scale=.8]{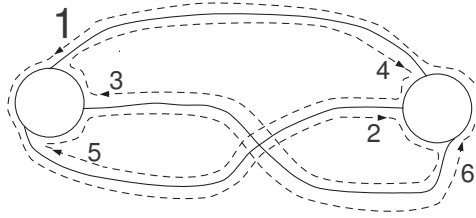}
	\caption{A fatgraph with $6$ half edges, the dashed curve represents its boundary component.}\label{1fig1}
\end{figure}

A fatgraph having $n$ edges and one boundary component can be encoded as a triple of permutations $(\alpha,\beta,\gamma)$
on $[2n]=\{1,2,\cdots 2n\}$ where $\gamma$ is the long cycle $(1,2,3,\ldots 2n)$.
This can be seen as follows: given a fatgraph $F$ with one boundary component,
we call the (two) ends of an edge half edges.
Pick a half edge, label it $1$, and start to travel the fatgraph counterclockwisely.
Label all visited half edges entering a vertex sequentially by labels $2,3,\ldots 2n$.
This induces two permutations $\alpha$ and $\beta$, where $\alpha$ is an involution without fixed points such that each
$\alpha$-cycle consists of the labels of the two half edges of the same edge and each $\beta$-cycle represents the
counterclockwise cyclic arrangement of all half edges incident to the same vertex.
By construction, $\gamma=(1,2,\ldots 2n)=\alpha\beta$, represents the unique boundary component of the fatgraph $F$.
An example of a fatgraph is illustrated in Figure~\ref{1fig1}, its corresponding triple of permutations are:
$$
\alpha=(1,4)(2,5)(3,6), \quad \beta=(1,5,3)(4,2,6), \quad \gamma=(1,2,3,4,5,6).
$$

A rooted one-face map is a one-face map where one half-edge is particularly marked and called the root.
We shall always label the root of a rooted one-face map with the label $1$, that is, we start from 
the root when traveling the boundary of the map.
Now given two rooted one-face maps which are respectively encoded into the triples
$(\alpha,\beta,\gamma)$ and $(\alpha',\beta',\gamma')$, they will be viewed as equivalent
if there exists a permutation $\pi$ such that
\begin{align*}
\alpha'=\pi\alpha \pi^{-1},\quad \gamma'=\pi\gamma\pi^{-1}, \quad \pi(1)=1.
\end{align*}
Following Euler's characteristic formula, the number of edges $n$, the number of vertices $v$ and
the genus $g$ of a one-face map satisfy
\begin{align*}
v-n+1=2-2g.
\end{align*}

The enumeration of rooted one-face maps has been extensively studied,
see for instance~\cite{chap2, chap1,chr-1,edmonds,gn,gs,hz,jac,wl} and the references therein.
Let $A(n,g)$ denote the number of rooted one-face maps (up to equivalence) of genus $g$ having $n$
edges and let $A_n(x)=\sum_{g\geq 0} A(n,g)x^{n+1-2g}$ be the corresponding generating function.
Four decades ago, Walsh and Lehman~\cite[eq.~$(13)$]{wl}, using a direct recursive method and
formal power series, obtained an explicit formula for $A(n,g)$ which can be reformulated as 
follows:
\begin{align}\label{1e2}
A(n,g)=\sum_{\lambda\vdash g}\frac{(n+1)n\cdots (n+2-2g-\ell(\lambda))}{2^{2g}
	\prod_i c_i!(2i+1)^{c_i}}\frac{(2n)!}{(n+1)!n!},
\end{align}
where the summation is taken over partitions $\lambda$ of $g$, $c_i$ is the number of parts $i$
in $\lambda$, and $\ell(\lambda)$ is the total number of parts.

More than a decade later, Harer and Zagier~\cite{hz} obtained in the context of computing the
virtual Euler characteristics of a curve:
\begin{align}\label{eq55}
A(n,g)=\frac{(2n)!}{(n+1)!(n-2g)!}[x^{2g}]\left(\frac{x/2}{\tanh x/2}\right)^{n+1},
\end{align}
where $[x^k]f(x)$ denotes the coefficient of $x^k$ in the expansion of the function $f(x)$.
Considering the relation between the RHS of eq.~\eqref{eq55} and its derivatives, they obtained
the following three-term recurrence, known as the Harer-Zagier recurrence:
\begin{align}\label{1eq4}
(n+1)A(n,g)  =2(2n-1)A(n-1,g)+(2n-1)(n-1)(2n-3)A(n-2,g-1).
\end{align}
They furthermore obtained the so-called Harer-Zagier formula:
\begin{align}\label{1eq5}
A_n(x)=\frac{(2n)!}{2^n n!}\sum_{k\geq 1}2^{k-1}{n\choose k-1}{x\choose k}.
\end{align}

There is a body of work on how to derive these results \cite{chap2,chap1,gn,gs,jac}.
A direct bijection for the Harer-Zagier formula was given in~\cite{gn}.
Combinatorial arguments to obtain the Lehman-Walsh formula and the Harer-Zagier recurrence
were recently given in~\cite{chap1}. One of the most recent advances is a new
recurrence for $A(n,g)$ obtained by Chapuy~\cite{chap2} via a bijective approach:
\begin{align}
2g A(n,g)= \sum_{k=1}^{g} {n+1-2(g-k) \choose 2k+1} A(n,g-k).
\end{align}
See also~\cite{chr-1} for a refinement of the recurrence and certain generalizations via plane permutations.

In this paper, we will prove
a new explicit formula for $A(n,g)$ and obtain a new formula for $A_n(x)$ which is more regular than the Harer-Zagier formula. We will also derive Chapuy's recursion.

A brief outline of the paper is as follows. In section~$2$, we shall employ an alternative interpretation of the Lehman-Walsh formula and construct two involutions
on pairs of permutations to obtain the following new formula for $A(n,g)$:
\begin{align}\label{1eq8}
A(n,g)=\frac{(2n)!}{2^n n!(n+1)!}\sum_{k=0}^n {n\choose k}\sum_{i+j=n+2-2g}C(n-k+1,i)(-1)^{k+1-j}C(k+1,j),
\end{align}
where $C(n,k)$ denotes the number of permutations on $n$ elements with $k$ cycles, i.e., the unsigned Stirling 
numbers of the first kind. This immediately gives us a new formula for the generating functions $A_n(x)$,~$n\geq 1$:
\begin{align}\label{1eq7}
A_n(x)=\frac{(2n)!}{2^n n!}\sum_{k\geq 0}{n\choose k}{x+n-k\choose n+1}.
\end{align}	
Utilizing the alternative interpretation of the Lehman-Walsh formula, another combinatorial explanation of the Harer-Zagier recurrence will be presented as well.

In section~$3$, by reformulating our expression 
for $A_n(x)$ in terms of the backward shift 
operator $\se: f(x)\rightarrow f(x-1)$ and proving a property satisfied by polynomials of the form
$p(\se)f(x)$, we easily establish Chapuy's recursion. Furthermore, by applying another property of polynomials of the form $p(\se)f(x)$ proved in Stanley~\cite{stan3}, we obtain the log-concavity of the numbers $A(n,g)$.



\section{New formulas for $A(n,g)$ and $A_n(x)$}

In the following, we will first prove a new formula for $A(n,g)$
by constructing two involutions on pairs of permutations.

We call a cycle of a permutation odd and even if it contains an odd and even number of elements, respectively.
Let $O(n+1,g)$ denote the number of permutations on $[n+1]$ which consist of $n+1-2g$ odd (disjoint) cycles.
For readers familiar with the formula for the number of permutations of a specific cycle type, see
Stanley~\cite[Prop.~$1.3.2$]{stan1}, it may be immediately realized that the Lehman-Walsh expression can be rewritten as
$$
\sum_{\lambda\vdash g}\frac{(n+1)n\cdots (n+2-2g-\ell(\lambda))}{2^{2g}\prod_i c_i!(2i+1)^{c_i}}\frac{(2n)!}{(n+1)!n!}=
\frac{(2n)!}{(n+1)!n!2^{2g}}O(n+1,g),
$$
thus we have more fundamental objects (permutations instead of maps) to work with.

Let $S=\{a, 1,2,\ldots n, b\}$. Let $\mathcal{T}_{A,l}$ denote the set of pairs
$(\alpha,\beta)$ where $\alpha$ is a permutation on $A\subset S$ while
$\beta$ is a permutation on $S\setminus A$, such that the sum of the number of $\alpha$- and $\beta$-cycles equals $l$.
Let
$$
\mathcal{T}_l=\bigcup_{A\subset S} \mathcal{T}_{A,l},
$$
where the union is taken over all $A\subset S$ such that $a \in A$ and $b\notin A$.
For each pair $(\alpha,\beta)\in \mathcal{T}_l$, we denote the difference between $\vert S\setminus A\vert$
and the number of $\beta$-cycles as $d(\beta)$ and set
$W[(\alpha,\beta)]=(-1)^{d(\beta)}$.
Then, it is clear that
\begin{lemma}\label{3lem7}
	\begin{align}
	\sum_{(\alpha,\beta)\in \mathcal{T}_l}W[(\alpha,\beta)]=
	\sum_{k=0}^n {n\choose k}\sum_{i+j=l}C(n-k+1,i)(-1)^{k+1-j}C(k+1,j),
	\end{align}
	where $C(n,k)$ denotes the number of permutations on $n$ elements with $k$ cycles.
\end{lemma}

Let $\mathcal{T'}\subset \mathcal{T}_l$ consist of pairs $(\alpha,\beta)$ where $\alpha(a)=a$ and $b$ is contained
in an odd cycle. We will construct our first involution $\phi$, which leads to

\begin{lemma}\label{3lem8}
	$$
	\sum_{(\alpha,\beta)\in \mathcal{T}_l} W[(\alpha,\beta)]= \sum_{(\alpha,\beta)\in \mathcal{T'}} W[(\alpha,\beta)].
	$$
\end{lemma}
\proof
Given $(\alpha,\beta)\in \mathcal{T}_l$, write both $\alpha$ and $\beta$ in their cycle decompositions and
denote the length of the cycle containing $b$ as $B$.
Define a map $\phi: \mathcal{T}_l\rightarrow \mathcal{T}_l$ as follows:
\begin{itemize}
	\item{Case~$1$:} if $(\alpha,\beta)\in \mathcal{T}'$, then $\phi: (\alpha,\beta)\mapsto (\alpha,\beta)$;
	\item{Case~$2$:} if $(\alpha,\beta)\notin \mathcal{T}'$, we distinguish two scenarios:
	\begin{itemize}
		\item If $B$ is odd and $\alpha(a)\neq a$, then $\phi:(\alpha,\beta)\mapsto (\alpha',\beta')$, where
		$\alpha'$ is obtained by deleting $\alpha(a)$ from the cycle decomposition of $\alpha$
		while $\beta'$ is obtained by inserting $\alpha(a)$ between $b$ and $\beta(b)$ in the cycle containing $b$
		and if $b=\beta(b)$, we map the cycle $(b)$ to $(b,\alpha(a))$.
		\item If $B$ is even then $b\neq \beta(b)$. Define $\phi: (\alpha,\beta)\mapsto (\alpha',\beta')$,
		where $\alpha'$ is obtained by inserting $\beta(b)$ between $a$ and $\alpha(a)$ and for $a=\alpha(a)$, we map
		$(a)$ to $(a, \beta(b))$. $\beta'$ is obtained by deleting $\beta(b)$ from the corresponding $\beta$-cycle.
	\end{itemize}
\end{itemize}	
We inspect that $\phi$ is an involution, whose fixed points are exactly all Case~$1$-pairs. Furthermore, $\phi$ preserves the
total number of cycles within the pairs. Accordingly, for Case~$1$-pairs, $\phi$ will preserve weights. For Case~$2$-pairs,
$\beta$ and $\beta'$ differ in the number of elements by $1$ but have the same number of cycles, so that $\phi$ will change
the sign of the weight, i.e., $W[(\alpha,\beta)]=-W[(\alpha',\beta')]$. Hence, the total weight of the $\phi$-orbit of any
Case~$2$-pair is $0$. Thus, $\sum_{(\alpha,\beta)\in \mathcal{T}_l} W[(\alpha,\beta)]$ is equal to the total weight of all
Case~$1$-pairs, completing the proof. \qed

Let $\mathcal{T''}$ be the set of pairs $(\alpha,\beta)\in \mathcal{T'}$ such that all cycles in $\alpha$ and $\beta$
are odd.
\begin{lemma}\label{3lem9}
	$$
	\sum_{(\alpha,\beta)\in \mathcal{T'}} W[(\alpha,\beta)]= \sum_{(\alpha,\beta)\in \mathcal{T''}} W[(\alpha,\beta)].
	$$
\end{lemma}


\proof
Define a map $\varphi: \mathcal{T'}\rightarrow \mathcal{T'}$ as follows:
\begin{itemize}
	\item{Case~$1$:} if $(\alpha,\beta)\in \mathcal{T}''$, then $\varphi: (\alpha,\beta) \mapsto (\alpha,\beta)$;
	\item{Case~$2$:} if $(\alpha,\beta)\notin \mathcal{T}''$,
	there is at least one even cycle in the collection of cycles from both $\alpha$ and $\beta$. Obviously, there
	is a unique even cycle, denoted by $C$, containing the
	minimal element among the union of elements from all even cycles.
	Let $\varphi: (\alpha, \beta) \mapsto (\alpha',\beta')$, where
	\begin{itemize}
		\item If $C$ is a cycle in $\alpha$, then $\alpha'=\alpha \setminus C$ and $\beta'=\beta \cup C$;
		\item Otherwise $\alpha'=\alpha \cup C$ and $\beta'=\beta \setminus C$.
	\end{itemize}
\end{itemize}
It is easy to see that $\varphi$ is an involution having Case-$1$ pairs as the only fixed points. Furthermore,
for Case-$1$ pairs, $\varphi$ preserves weights. For Case-$2$ pairs, $\beta$ and $\beta'$
differ in the number of elements by an even number while they differ in the number of cycles by $1$.
Hence $W[(\alpha,\beta)]=-W[\varphi((\alpha,\beta))]$.
As a result, the total weight over $(\alpha,\beta)$ in Case~$2$ is $0$, completing the proof. \qed

Based on these lemmas above, we obtain
\begin{theorem}
	For $n,l\geq 0$, we have
	\begin{align}\label{2e6}
	\sum_{k=0}^n {n\choose k}\sum_{i+j=l}C(n-k+1,i)(-1)^{k+1-j}C(k+1,j)=2^{l-2}O(n+1,\frac{n+2-l}{2}).
	\end{align}
\end{theorem}
\proof
Lemma~\ref{3lem7}, Lemma~\ref{3lem8} and Lemma~\ref{3lem9} imply, that the LHS of eq.~\eqref{2e6} equals
$\sum_{(\alpha,\beta)\in \mathcal{T''}} W[(\alpha,\beta)]$, where the total number of cycles in $\alpha$ and $\beta$ is $l$.
Since in $\mathcal{T''}$ all cycles are odd, the number of elements and the number of cycles in $\beta$ have the same
parity. Thus, for any $(\alpha,\beta)\in \mathcal{T''}$, $W[(\alpha,\beta)]=1$, i.e~the total weight over
$\mathcal{T''}$ equals the total number of elements in $\mathcal{T''}$.

Since $a$ is a fixed point in $\alpha$ for any $(\alpha, \beta)\in \mathcal{T''}$, each pair $(\alpha, \beta)\in
\mathcal{T''}$ can be viewed as a partition of all $l-1$ odd cycles of a permutation on
$[n]\cup \{b\}$, except the cycle containing $b$, into two ordered parts.
Conversely, given a permutation on $[n]\cup \{b\}$ with $l-1$ odd cycles, there are $2^{l-2}$
different ways to partition all cycles except the one containing $b$ into two ordered parts.
Therefore, we have
\begin{align*}
\sum_{(\alpha,\beta)\in \mathcal{T''}}W[(\alpha,\beta)]=|\mathcal{T''}|=2^{l-2}O(n+1,\frac{n+2-l}{2}),
\end{align*}
completing the proof. \qed

As a corollary, we obtain a new recurrence for the unsigned Stirling numbers of the first kind:
\begin{corollary}
	For $n\geq 0$ and $l\neq n \mod 2$, we have
	\begin{align}\label{2e7}
	\sum_{k=0}^n {n\choose k}\sum_{i+j=l}(-1)^{k+1-j}C(n-k+1,i)C(k+1,j)=0.
	\end{align}
\end{corollary}

\proof
Since every cycle is odd for any $(\alpha, \beta)\in \mathcal{T}''$, the number of total elements $n+2$
has the same parity as the total number of cycles $l$. If $l\neq n \mod 2$, $\mathcal{T}''=\emptyset$, i.e., $O(n+1,\frac{n+2-l}{2})=0$, whence
\begin{align*}
\sum_{k=0}^n {n\choose k}\sum_{i+j=l}(-1)^{k+1-j}C(n-k+1,i)C(k+1,j)=0.
\end{align*}
\qed

Setting $l=n+2-2g$ in eq.~\eqref{2e6}, we have
\begin{corollary}
	For $n,g\geq 0$, we have
	\begin{align}\label{2e8}
	\sum_{k=0}^n {n\choose k}\sum_{i+j=n+2-2g}C(n-k+1,i)(-1)^{k+1-j}C(k+1,j)=2^{n-2g}O(n+1,g).
	\end{align}
\end{corollary}

Accordingly, multiplying $\frac{(2n)!}{2^n n!(n+1)!}$ on both sides of eq.~\eqref{2e8} we obtain a new
explicit formula for $A(n,g)$:

\begin{theorem}\label{2thm2}
	$A(n,g)=\frac{(2n)!}{2^n n!(n+1)!}\bar{A}(n,g)$ where
	\begin{align}\label{3e15}
	\bar{A}(n,g)=\sum_{k=0}^n {n\choose k}\sum_{i+j=n+2-2g}C(n-k+1,i)(-1)^{k+1-j}C(k+1,j).
	\end{align}
\end{theorem}


Then, we immediately obtain

\begin{theorem}\label{L:11} The generating functions $A_n(x)$ for $n\geq 0$ satisfy
	\begin{align}\label{2e4}
	\sum_{g\geq 0} A(n,g)x^{n+1-2g}=\frac{(2n)!}{2^n n!}\sum_{k\geq 0}{n\choose k}{x+n-k\choose n+1}.
	\end{align}
\end{theorem}
\proof
According to \cite{stan1} we have
\begin{align*}
x(x+1)(x+2)\cdots (x+n-1)&=\sum_{k\geq 1} C(n,k)x^k,\\
x(x-1)(x-2)\cdots (x-n+1)&=\sum_{k\geq 1} (-1)^{n-k}C(n,k)x^k.
\end{align*}
From these facts and eq.~\eqref{3e15},
we immediately compute
\begin{align*}
&\bar{A}(n,g)\nonumber\\
&=\sum_{k=0}^n {n\choose k}[x^{n+2-2g}][(x+n-k)(x+n-k-1)\cdots x][x(x-1)\cdots (x-k)]\\
&=\sum_{k=0}^n {n\choose k}[x^{n+1-2g}](x+n-k)(x+n-k-1)\cdots x(x-1)\cdots (x-k),
\end{align*}
from which the theorem follows. \qed

\begin{remark}
	Clearly, depending on $n$, $A_n(x)$ represents, either an odd or even function, which is not obvious
	from Harer-Zagier's formula.
	Our new formula on the RHS of eq.~\eqref{2e4} makes this feature immediately evident:
	let $(a)_n$ denote the falling
	factorial $a(a-1)\cdots (a -n +1)$.
	Then,
	\begin{align*}
	{n\choose k}\frac{(x+n-k)_{n+1}}{(n+1)!}&={n\choose k} \frac{(x+n-k)(x+n-k-1)\cdots (x-k)}{(n+1)!}\\
	&=(-1)^{n+1}{n\choose k} \frac{[-(x+n-k)][-(x+n-k-1)]\cdots [-(x-k)]}{(n+1)!}\\
	&=(-1)^{n+1}{n\choose n-k} \frac{(-x+n-(n-k))_{n+1}}{(n+1)!},
	\end{align*}
	which implies $A_n(x)=(-1)^{n+1}A_n(-x)$.
\end{remark}

For the three-term Harer-Zagier recurrence eq.~\eqref{1eq4}, the proof in~\cite{chap1} is based on 
a combinatorial isomorphism. Here, we give another indirect combinatorial proof for it, which relies on the three-term recurrence for $O(n,g)$ evident in the following lemma.

\begin{lemma}
	The numbers $O(n,g)$ satisfy
	\begin{align}\label{a2}
	O(n+1,g)=O(n,g)+n(n-1)O(n-1,g-1).
	\end{align}
\end{lemma}
\proof
Firstly, we claim that

\emph{Claim.} The numbers $O(n,g)$ satisfy
\begin{align}\label{a1}
O(n+1,g)=O(n,g)+\sum_{k=1}^g (n)_{2k}O(n-2k,g-k).
\end{align}
This can be seen by classifying all permutations on $[n+1]$ with $n+1-2g$ odd
cycles based on the length $\ell_1$ of the cycle containing the element $1$.
Obviously, the class with $\ell_1=1$ has $O(n,g)$ permutations.
For the class with $\ell_1=2k+1$, $k\geq 1$, there are $(n)_{2k}$ ways to choose $2k$
elements from $[n+1]\setminus \{1\}$ which together with the element $1$ form a cycle of length $2k+1$. The remaining $n-2k$ elements
can arbitrarily form $n-2g=(n-2k)-2(g-k)$ odd cycles, and there are $O(n-2k,g-k)$ different
ways to do that. Thus, eq.~\eqref{a1} follows.

Next, we note eq.~\eqref{a1}  implies
\begin{align*}
O(n-1,g-1)=O(n-2,g-1)+\sum_{k=1}^{g-1} (n-2)_{2k}O(n-2k-2,g-1-k).
\end{align*}
Thus we obtain
\begin{align*}
O(n+1,g)&=O(n,g)+\sum_{k=1}^g (n)_{2k}O(n-2k,g-k)\nonumber\\
&=O(n,g)+n(n-1)[O(n-2,g-1)+\sum_{k=2}^g (n-2)_{2k-2}O(n-2k,g-k)]\nonumber\\
&=O(n,g)+n(n-1)O(n-1,g-1),
\end{align*}
completing the proof of the lemma. \qed

\begin{proposition}[Harer-Zagier recurrence]
	\begin{align*}
	(n+1)A(n,g)  =2(2n-1)A(n-1,g)+(2n-1)(n-1)(2n-3)A(n-2,g-1).
	\end{align*}
\end{proposition}
\proof
Using eq.~\eqref{a2}, we have
\begin{align*}
A(n,g)&=\frac{(2n)!}{(n+1)!n!2^{2g}}O(n+1,g)\\
&=\frac{(2n)!}{(n+1)!n!2^{2g}}O(n,g)+\frac{(2n)!}{(n+1)!n!2^{2g}}n(n-1)O(n-1,g-1)\\
&= \frac{2n(2n-1)}{(n+1)n}\cdot \frac{(2n-2)!}{n!(n-1)!2^{2g}}O(n,g)\\
&\quad +\frac{2n(2n-1)(2n-2)(2n-3)}{(n+1)nn(n-1)2^2}\cdot \frac{(2n-4)!n(n-1)}{(n-1)!(n-2)!2^{2g-2}}O(n-1,g-1)\\
&= \frac{2(2n-1)}{n+1}A(n-1,g)+\frac{(2n-1)(n-1)(2n-3)}{n+1}A(n-2,g-1),
\end{align*}
whence the proposition. \qed

\section{Chapuy's recursion and log-concavity}

In this section, by reformulating our expression 
for $A_n(x)$ in terms of the backward shift 
operator $\se: f(x)\rightarrow f(x-1)$ and proving a property satisfied by polynomials of the form
$p(\se)f(x)$, we easily establish Chapuy's recursion. Furthermore, by applying another property of polynomials of the form $p(\se)f(x)$ proved in Stanley~\cite{stan3}, we obtain the log-concavity of the numbers $A(n,g)$.

First, our new formula of $A_n(x)$ implies

\begin{proposition}\label{P:13}
	\begin{align}\label{2e66}
	\sum_{g\geq 0} A(n,g)x^{n+1-2g}=\frac{(2n)!}{2^n n! (n+1)!}(1+\mathrm{E} )^n (x+n)_{n+1}.
	\end{align}
\end{proposition}
\proof This is evident from the following computation:
\begin{align*}
\sum_{g\geq 0} A(n,g)x^{n+1-2g}&=\frac{(2n)!}{2^n n!}\sum_{k\geq 0}{n\choose k}{x+n-k\choose n+1}\\
&=\frac{(2n)!}{2^n n!(n+1)!}\sum_{k\geq 0}{n\choose k}(x-k+n)_{n+1}\\
&=\frac{(2n)!}{2^n n!(n+1)!}\sum_{k\geq 0}{n\choose k}\se^k(x+n)_{n+1}\\
&=\frac{(2n)!}{2^n n!(n+1)!}(1+{\se})^n (x+n)_{n+1}.
\end{align*}\qed

We proceed by showing that any polynomial of the form $p(\se)(x+n)_{n+1}$ satisfies:
\begin{theorem}\label{4tgeneral}
	Let $p(t)=\sum_{k=0}^n a_k t^k$ and $F(x)=p(\mathrm{E})(x+n)_{n+1}$.
	If $\frac{a_1}{a_0}=\frac{a_{n-1}}{a_n}$ and
	$
	k a_k +(n-k+2) a_{k-2}=\frac{a_1}{a_0} a_{k-1}$, {for $2\leq k \leq n$},
	then
	\begin{align}\label{general1}
	(n+2+\frac{a_1}{a_0})F(x)=x(F(x+1)-F(x-1)).
	\end{align}
	Moreover, let $b_k= [x^k] p(\mathrm{E})(x+n)_{n+1}$, then we have
	\begin{align}\label{general2}
	\left( \frac{n+2+\frac{a_1}{a_0}}{2}-k\right) b_k=\sum_{j\geq 1} {k+2j \choose 2j+1}b_{k+2j}.
	\end{align}
\end{theorem}

\proof
Note that by assumption, the RHS of eq.~\eqref{general1} is equal to
$$
\sum_{k=0}^n \{a_k x {(x+1+n-k)_{ n+1}} - a_k x {(x-1+n-k)_{ n+1}}\}.
$$
Clearly, we have
\begin{align*}
a_k x {(x+1+n-k)_{ n+1}}&=a_k (x-k+k) {(x+1+n-k)_{ n+1}}\\
&=(x+1+n-k)a_k {(x+n-k)_{ n+1}}+ k a_k{(x+1+n-k)_{ n+1}},\\
a_k x {(x-1+n-k)_{ n+1}}&=a_k[x+(n-k)-(n-k)]{(x-1+n-k)_{ n+1}}\\
&=(x-k-1)a_k{(x+n-k)_{ n+1}}-(n-k)a_k {(x-1+n-k)_{ n+1}}.
\end{align*}
Then, to obtain eq.~\eqref{general1}, it suffices to show that the difference between the respective sums of the RHS of the last two equations equals the LHS of eq.~\eqref{general1}. This follows from the following computations:
\begin{align*}
& \sum_{k=0}^n (x+1+n-k)a_k {(x+n-k)_{ n+1}} - \sum_{k=0}^n (x-k-1)a_k{(x+n-k)_{ n+1}}\\
=& (n+2) \sum_{k=0}^n a_k {(x+n-k)_{ n+1}}= (n+2) F(x),
\end{align*}
and
\begin{align*}
& \sum_{k=0}^n k a_k{(x+1+n-k)_{ n+1}} + \sum_{k=0}^n (n-k)a_k {(x-1+n-k)_{ n+1}}\\
=& \sum_{k=2}^n \left\{k a_k{(x+1+n-k)_{ n+1}}+  [n-(k-2)]a_{k-2} {(x-1+n-(k-2))_{ n+1}} \right\}+\\
&\quad\quad\quad  a_1{(x+1+n-1)_{n+1}}+ [n-(n-1)]a_{n-1}{(x-1+n-(n-1))_{n+1}}\\
=& \sum_{k=2}^n \left[ k a_k +(n-k+2) a_{k-2}\right]{(x+n-(k-1))_{n+1}}+ a_1{(x+n)_{n+1}}+ a_{n-1}{(x)_{n+1}}\\
=& \frac{a_1}{a_0} \sum_{k=1}^{n-1} a_k {(x+n-k)_{n+1}}+ \frac{a_1}{a_0} a_0{(x+n)_{n+1}}+\frac{a_1}{a_0} a_n{(x)_{n+1}} =\frac{a_1}{a_0} F(x).
\end{align*}
Next, the polynomial $F(x)$ is analytic and has thus a power series expansion everywhere. In particular, we have
\begin{align*}
F(x+1) =\sum_ {k\geq 0} \frac{F^{(k)}(x)}{k!}(x+1-x)^k, \quad
F(x-1) =\sum_ {k\geq 0} \frac{F^{(k)}(x)}{k!}(x-1-x)^k.
\end{align*}
Then,
$$
\frac{n+2+\frac{a_1}{a_0}}{2}F(x)=\frac{x (F(x+1)-F(x-1))}{2}=\sum_{k\geq 0} \frac{x F^{(2k+1)}(x)}{(2k+1)!},
$$
which can be reformulated as
$$
\frac{n+2+\frac{a_1}{a_0}}{2}F(x)-xF'(x)=\sum_{k\geq 0} \left(	\frac{n+2+\frac{a_1}{a_0}}{2}-k\right) b_k x^{k} =\sum_{j\geq 1} \frac{x F^{(2j+1)}(x)}{(2j+1)!}.
$$
Comparing the coefficients of the last equation based on the fact that
$$
\frac{xF^{(2j+1)}(x)}{(2j+1)!}=\sum_{i\geq 0} \frac{(i)_{2j+1}}{(2j+1)!} b_ix^{i-2j},
$$
we obtain eq.~\eqref{general2} and the proof of the theorem is complete.\qed

As a consequence of Theorem~\ref{4tgeneral}, we immediately obtain Chapuy's recurrence.
\begin{corollary}[Chapuy's recursion]
	\begin{align}\label{4eq11}
	2g A(n,g) =\sum_{k=1}^{g} {n+1-2(g-k) \choose 2k+1} A(n, g-k).
	\end{align}
\end{corollary}
\proof
For $A_n(x)$, Proposition~\ref{P:13} gives us $p(t)=\sum_{k=0}^n a_k t^k$ where $a_k=\frac{(2n)!}{2^n n! (n+1)!}{n \choose k}$. It is obvious that $\frac{a_1}{a_0}=\frac{a_{n-1}}{a_n}=n$. Furthermore,
\begin{align*}
k a_k +(n-k+2) a_{k-2} &=\frac{(2n)!}{2^n n! (n+1)!} \left[ k{n\choose k}+  [n-(k-2)]{n\choose k-2} \right]\\
&=\frac{(2n)!}{2^n n! (n+1)!} \left[ n {n-1 \choose k-1}+ n {n-1 \choose k-2} \right]\\
&=\frac{(2n)!}{2^n n! (n+1)!} n{n \choose k-1}=na_{k-1}.
\end{align*}
Hence, we can apply Theorem~\ref{4tgeneral} to $A_n(x)$ and obtain 
\begin{align*}
2g A(n,g)&=\left(\frac{n+2+n}{2}-(n+1-2g)\right)[x^{n+1-2g}]A_n(x)\\
&=\sum_{j\geq 1} {n+1-2g+2j \choose 2j+1}[x^{n+1-2g+2j}] A_n(x)\\
&=\sum_{k=1}^{g} {n+1-2(g-k) \choose 2k+1} A(n, g-k),
\end{align*}
which is Chapuy's recurrence. \qed

\subsection{Log-concavity of $A(n,g)$}

Proposition~\ref{P:13} has a close relative in \cite[eq.~$(6)$]{stan3} in terms of a formula using the backward shift operator
which can be reformulated as
\begin{align}\label{2eq77}
\sum_{g\geq 0} A(n,g)x^{n+1-2g}=\frac{1}{2^n n! (2n+1)}(1-\se^2)^n (x+2n)_{2n+1}.
\end{align}
This formula was obtained using the character theory of the symmetric group.
In addition, it was proved in~\cite{stan3} that the RHS of the last equation is a degree $n+1$ polynomial 
despite it appears to be a degree $2n+1$ polynomial.

In~\cite{stan3}, it is further proved that every zero of the polynomial on the RHS of eq.~\eqref{2eq77}
is purely imaginary, based on the following result:
\begin{proposition}[Stanley~\cite{stan3}]\label{2prop111}
	Let $p(t)$ be a complex polynomial of degree exactly $d$, such
	that every zero of $p(t)$ lies on the circle $|z| = 1$. Suppose that the multiplicity
	of $1$ as a root of $p(t)$ is $m \geq 0$. Let $P(x) = p(\mathrm{E})(x + n -1)_n$.
	If $d\leq n-1$, then
	$$
	P(x)=(x+n-d-1)_{n-d}Q(x),
	$$
	where $Q(x)$ is a polynomial of degree $d-m$ for which every zero has
	real part $\frac{d-n+1}{2}$.
\end{proposition}

Applying Proposition~\ref{2prop111} to the RHS of eq.~\eqref{2e66}, we have $p(t)=(1+t)^n$,~$d=n$,~$m=0$ so that $\sum_{g\geq 0} A(n,g)x^{n+1-2g}=xQ(x)$ with
every zero of $Q(x)$ being purely imaginary. Let $H(x)=\sum_{g\geq 0} A(n,g)x^{\lfloor \frac{n+1}{2}\rfloor-g}$.
Then,
we have

\begin{corollary}
	$H(x)$ has only non-positive real zeros, and for fixed $n\geq 1$, the sequence $A(n,g)$ is log-concave, i.e., $A(n,g)^2\geq
	A(n,g-1)A(n,g+1)$.
\end{corollary}

\proof
By construction, we have
either $H(x^2)=xQ(x)$ or $H(x^2)=Q(x)$, depending on the parity of $n$.
In either case, it implies $H(x)$ has only non-positive real zeros.
It is well known that having only real zeros implies log-concavity of the coefficients (see~\cite{stan4} for instance). This completes the proof. \qed





\end{document}